\documentclass[12pt,english]{article}
\usepackage[T1]{fontenc}
\usepackage[latin9]{inputenc}
\usepackage{geometry}
\usepackage{amsmath}
\usepackage{graphicx}
\usepackage{setspace}
\usepackage[authoryear]{natbib}
\makeatletter

\providecommand{\tabularnewline}{\\}

\newcommand{\lyxaddress}[1]{
	\par {\raggedright #1
	\vspace{1.4em}
	\noindent\par}
}

\@ifundefined{date}{}{\date{}}
\usepackage{pdflscape}

\makeatother

\usepackage{babel}
\begin{document}
\title{Convergent estimators of variance of a spatial mean in the presence
of missing observations}
\author{Ashwin K Seshadri}
\maketitle

\lyxaddress{Centre for Atmospheric and Oceanic Sciences and Divecha Centre for
Climate Change, Indian Institute of Science, Bangalore 560012, India,
email: ashwins@iisc.ac.in}
\begin{abstract}
In the geosciences, a recurring problem is one of estimating spatial
means of a physical field using weighted averages of point observations.
An important variant is when individual observations are counted with
some probability less than one. This can occur in different contexts:
from missing data to estimating the statistics across subsamples.
In such situations, the spatial mean is a ratio of random variables,
whose statistics involve approximate estimators derived through series
expansion. The present paper considers truncated estimators of variance
of the spatial mean and their general structure in the presence of
missing data. To all orders, the variance estimator depends only on
the first and second moments of the underlying field, and convergence
requires these moments to be finite. Furthermore, convergence occurs
if either the probability of counting individual observations is larger
than $1/2$ or the number of point observations is large. In case
the point observations are weighted uniformly, the estimators are
easily found using combinatorics and involve Stirling numbers of the
second kind. 
\end{abstract}

\section{Introduction}

Very often in the geosciences, we seek estimates of the spatial average
of some physical variable such as rainfall or temperature (\citet{Rodriguez-Iturbe1974,Shen1994,Morrissey1999,Kundu2007,Villarini2008,Prakash2019}).
Examining the properties of such an average is a decades-old problem
(\citet{Kagan1997}). An interesting variant is where individual observations
are reported with known probability that is less than one. In this
case, what form do the estimators of bias and variance of the spatial
average take? \citet{Seshadri2018} examined such a situation, involving
point observations $r_{i}$ at locations $i$ $\left(1\leq i\leq N\right)$,
with each observation assigned weight $\beta_{i}$, such that $\sum_{i}\beta_{i}=1$.
Time is assumed discrete, and not all observations are reported at
each time. The reporting status is measured by Boolean random variable
$s_{i}$, with $s_{i}=1$ if reported at that time and $s_{i}=0$
otherwise. In our setting, the reporting status of each observation
site has been assumed to be independent of whether any other sites,
including those in close proximity, are reported. This makes the $s_{i}$'s
statistically independent. Assuming that each $s_{i}$ has probability
$\alpha$ of equaling one, we define the spatial average over samples
without replacement
\begin{equation}
r\left(t\right)\equiv\frac{R\left(t\right)}{S\left(t\right)}=\frac{\sum_{i}\beta_{i}s_{i}\left(t\right)r_{i}\left(t\right)}{\sum_{i}\beta_{i}s_{i}\left(t\right)},\label{eq:p1}
\end{equation}
which describes the weighted mean over only those observations that
are reported. The weights $\beta_{i}$ are constant in time, but $s_{i}$
and $r_{i}$ are time-varying. Henceforth, we omit the time-dependence
from the formulas for clarity of expression, and denote the spatial
average as $r=f\left(R,S\right)=R/S$. In general, whenever $f\left(R,S\right)$
is analytic in some neighborhood of $f\left(\boldsymbol{\mathbf{E}}R,\boldsymbol{\mathbf{E}}S\right),$
with $\boldsymbol{\mathbf{E}}$ denoting expectation, a Taylor series
expansion can be made for $f$ in the neighborhood and is convergent. 

The above ratio describes the spatial mean over available observations.
Different realizations of the set $\left\{ s_{i}:i=1,2,\ldots N\right\} $
give rise to different estimates of the mean, and statistics across
the possible realizations are often of interest. It is central to
this problem that $s_{i}$ is a random variable and origins of randomness
in counting individual observations can be manifold. In general such
randomness of $s_{i}$ invokes two types of effects. Individual observations
might be missing or unreported with nonzero probability, in which
case the spatial average at each time involves only those observations
reported at that time (\citet{Vinnikov2004}), and one might be interested
in statistics of the reported spatial average, such as its temporal
variance (\citet{Seshadri2018}). Alternately, the random variable
$s_{i}$ might invoke deliberate omissions of observations to estimate
properties of sub-samples without replacement (\citet{Prakash2019}).
In the first case, the variance estimator would count additional effects
of missing data on the variance. In the second case, it would provide
a nonparametric estimate of the dispersion of the ensemble that is
comprised of different sub-samples. Formally, however, both problems
are the same. 

Therefore, while the formula in Eq. (\ref{eq:p1}) is simple, perhaps
deceptively so, its statistics recur in many guises. Moreover, if
the $s_{i}$'s are fixed to be equal to one, standard formulas for
bias and variance reemerge (\citet{Gandin1993,Kagan1997}). However,
the interesting situation arises from randomness in these variables,
where the denominator is also a random variable so that there is no
exact formula for the statistics of Eq. (\ref{eq:p1}). As shown in
\citet{Seshadri2018}, and following a much earlier literature that
is broader in scope (\citet{Hartley1954,Goodman1958,Hinkley1969,Kendall1994}),
useful approximations can nonetheless been derived. Such approximate
formulas appear to perform remarkably well when benchmarked against
Monte Carlo simulations, despite being based on low-order truncations
(\citet{Seshadri2018,Prakash2019}). Therefore the convergence of
these estimators derived from power series (\citet{Oehlert1992})
merits inquiry, and such is the goal of the present paper. 

Our basic approach is conceptually straightforward and involves taking
expectations of a truncated Taylor expansion of the function (\citet{Oehlert1992}).
Successive degrees of approximation gives rise to corresponding estimators.
The main question about this approach is whether such estimators converge,
and under what conditions they do so, as one progressively includes
additional terms in the Taylor series. The basic theory is well-known
and there are two requirements. First, the Taylor series itself must
converge (\citet{Bromwich1991}). Its convergence might be conditional,
valid only for small distances from the expansion center. Or it might
be unconditional, valid for all values of the underlying field. Either
situations are possible in the more general problem of deriving approximate
estimators by this approach, whose examples range beyond the spatial
geosciences to various applied fields in statistics and economics
(\citet{Loistl1976,Hlawitschka1994,Markowitz2015}). Generally, the
Taylor series converges only for a finite radius (\citet{Bromwich1991,Gemignani1970,Hlawitschka1994}).
Within this radius of convergence, the further question is whether
the expectation also converges. This is governed by the dominated
convergence theorem (\citet{Weir1973}). Conditions for dominated
convergence are met if the Taylor series converges absolutely, in
which case the expectation is also convergent. This will be the relevant
test for our estimators, which involve the expectation of series that
describe the variance of the spatial mean in Eq. (\ref{eq:p1}). 

\section{Theory and derivations}

\subsection{Convergence}

Following \citet{Seshadri2018}, we consider point observations $r_{i}$
of a spatially varying field at locations indexed by $i$ $\left(1\leq i\leq N\right)$.
Individual observations are missing with constant probability $1-\alpha$,
with $0<\alpha\leq1$. The status of each observation is described
by Boolean random variable $s_{i}$: if the observation is reported
then $s_{i}=1$ otherwise $s_{i}=0$. The spatial average is defined
in Eq. (\ref{eq:p1}) and this paper seeks a formula for its variance.
Since both numerator and denominator are random variables, there is
no exact expression and we must expand $f\left(R,S\right)\equiv R/S$
by its Taylor series about $\left(\boldsymbol{\mathbf{E}}R,\boldsymbol{\mathbf{E}}S\right)$,
with $\boldsymbol{\mathbf{E}}$ denoting expectation
\begin{multline}
f\left(R,S\right)=f\left(\boldsymbol{\mathbf{E}}R,\boldsymbol{\mathbf{E}}S\right)+\sum_{n=1}^{\infty}\frac{1}{n!}\sum_{k=0}^{n}\left(\begin{array}{c}
n\\
k
\end{array}\right)\frac{\partial^{n}f}{\partial R^{n-k}\partial S^{k}}\left(R-\boldsymbol{\mathbf{E}}R\right)^{n-k}\left(S-\boldsymbol{\mathbf{E}}S\right)^{k}\label{eq:p2}
\end{multline}
where $\left(\begin{array}{c}
n\\
k
\end{array}\right)=\frac{n!}{\left(n-k\right)!k!}$. This gives rise to series expansion
\begin{equation}
f\left(R,S\right)=\sum_{n=0}^{\infty}\left(-1\right)^{n}\left(\frac{1}{\boldsymbol{\mathbf{E}}S}\right)^{n+1}\left(R-\boldsymbol{\mathbf{E}}R\right)\left(S-\boldsymbol{\mathbf{E}}S\right)^{n}+\sum_{n=0}^{\infty}\left(-1\right)^{n}\frac{\boldsymbol{\mathbf{E}}R}{\left(\boldsymbol{\mathbf{E}}S\right)^{n+1}}\left(S-\boldsymbol{\mathbf{E}}S\right)^{n}\label{eq:p3}
\end{equation}
as shown in Supplementary Information (SI) section 1. Writing this
series as
\begin{equation}
f\left(R,S\right)=\sum_{n=0}^{\infty}a_{n}+\sum_{n=0}^{\infty}b_{n}\equiv\sum_{n=0}^{\infty}c_{n},\label{eq:p4}
\end{equation}
its square is the Cauchy product of the series with itself
\begin{equation}
f^{2}\left(R,S\right)=\left(\sum_{n=0}^{\infty}c_{n}\right)\left(\sum_{n=0}^{\infty}c_{n}\right).\label{eq:p5}
\end{equation}
When does the above series for $f^{2}\left(R,S\right)$ converge?
According to Merten's theorem, a Cauchy product converges if the individual
series are convergent and at least one of the series converges absolutely
(\citet{Rudin1976}). Therfore, it is sufficient that $\sum_{n=0}^{\infty}c_{n}$
converges absolutely. From the triangle inequality 
\begin{equation}
\sum_{n=0}^{\infty}\left|c_{n}\right|=\sum_{n=0}^{\infty}\left|a_{n}+b_{n}\right|\leq\sum_{n=0}^{\infty}\left|a_{n}\right|+\sum_{n=0}^{\infty}\left|b_{n}\right|\label{eq:p6}
\end{equation}
it is sufficient, for absolute convergence of $\sum_{n=0}^{\infty}c_{n}$,
that $\sum_{n=0}^{\infty}a_{n}$ and $\sum_{n=0}^{\infty}b_{n}$ converge
absolutely. Let us evaluate the absolute convergence of each of the
two series. For the first one
\begin{equation}
\sum_{n=0}^{\infty}\left|a_{n}\right|=\left|R-\boldsymbol{\mathbf{E}}R\right|\sum_{n=0}^{\infty}\left(\frac{1}{\boldsymbol{\mathbf{E}}S}\right)^{n+1}\left|S-\boldsymbol{\mathbf{E}}S\right|^{n},\label{eq:p7}
\end{equation}
where we have used the fact that $\boldsymbol{\mathbf{E}}S>0$ since
$\alpha>0$. A sufficient condition for its convergence is provided
by the ratio test. We require the ratio 
\begin{equation}
\frac{\left|a_{n+1}\right|}{\left|a_{n}\right|}=\frac{\left|S-\boldsymbol{\mathbf{E}}S\right|}{\boldsymbol{\mathbf{E}}S}<1\label{eq:p8}
\end{equation}
as $n\rightarrow\infty$. Since the weights $\beta_{i}$ in Eq. (\ref{eq:p1})
sum to one, the expectation is $\boldsymbol{\mathbf{E}}S=\sum_{i=1}^{N}\beta_{i}\boldsymbol{\mathbf{E}}s_{i}=\sum_{i=1}^{N}\beta_{i}\alpha=\alpha$.
Moreover, the random variable $S$ is bounded by $0<S\leq1$, with
$S=0$ precluded because this would entail the absence of any reported
observations. This bound ensures that the above ratio condition is
met if
\begin{equation}
\frac{1-\alpha}{\alpha}<1,\label{eq:p9}
\end{equation}
or $\alpha>1/2$. 

In case $\alpha<1/2$, convergence is not assured in case $S\geq2\boldsymbol{\mathbf{E}}S$.
The probability of this occurrence can be bounded through Hoeffding's
inequality. We use a version of this result for bounded random variables
that are independent but not identically distributed. We can write
$S=\sum_{i}x_{i}$, where $x_{i}=\beta_{i}s_{i}$. The individual
$x_{i}$'s are bounded between $m_{i}$ and $M_{i}$, and are independent,
since the $s_{i}$'s are independent. Hoeffding's result shows that
the probability of the event that $\sum_{i}\left(x_{i}-\boldsymbol{\mathbf{E}}x_{i}\right)\geq c$
is less than or equal to $e^{-2c^{2}/\sum_{i}\left(M_{i}-m_{i}\right)^{2}}$
(\citet{Vershynin2018}). Since $0\leq s_{i}\leq1$, $m_{i}=0$ and
$M_{i}=\beta_{i}$. Moreover $\boldsymbol{\mathbf{E}}x_{i}=\alpha\beta_{i}$,
so $\sum_{i}\boldsymbol{\mathbf{E}}x_{i}=\alpha\sum_{i}\beta_{i}=\alpha$.
This result implies that the probability of $S-\alpha\geq\alpha$
, or equivalently $S\geq2\boldsymbol{\mathbf{E}}S$, is smaller than
\begin{equation}
e^{-2\alpha^{2}/\sum_{i}\beta_{i}^{2}},\label{eq:p10}
\end{equation}

with stricter bounds arising if the distribution of weights is more
uniform, making $\sum_{i}\beta_{i}^{2}$ smaller. In the limiting
case of uniform weights $\beta_{i}=1/N$ and $\sum_{i=1}^{N}\beta_{i}^{2}=1/N$,
the probability that $S>2\boldsymbol{\mathbf{E}}S$ is bounded above
by $e^{-2N\alpha^{2}}$. In practice the Hoeffding inequality is quite
conservative, and stronger bounds are possible, especially for small
$\alpha$, but the calculation does illustrate how convergence is
favored by large $N$. 

For uniform weights, bounds can be found more directly, since $S$
takes a binomial distribution with mean $\alpha$ and standard deviation
$\sqrt{\frac{\alpha\left(1-\alpha\right)}{N}}$. Therefore the value
that $S$ takes must be $\sqrt{N\frac{\alpha}{1-\alpha}}$ standard
deviations away from the mean, which is very unlikely for large $N$.
For example, even with a small probability $\alpha=0.1$, if $N=100$
then $S$ must be $3.3$ standard deviations away. In practice, for
geophysical problems the probability of reporting observations is
considerably larger. 

Similarly, for absolute convergence of the second series in Eq. (\ref{eq:p6})
\begin{equation}
\sum_{n=0}^{\infty}\left|b_{n}\right|=\left|\boldsymbol{\mathbf{E}}R\right|\sum_{n=0}^{\infty}\frac{1}{\left(\boldsymbol{\mathbf{E}}S\right)^{n+1}}\left|S-\boldsymbol{\mathbf{E}}S\right|^{n},\label{eq:p11}
\end{equation}

the condition 
\begin{equation}
\frac{\left|b_{n+1}\right|}{\left|b_{n}\right|}=\frac{\left|S-\boldsymbol{\mathbf{E}}S\right|}{\boldsymbol{\mathbf{E}}S}<1\label{eq:p12}
\end{equation}
is identical. In summary, absolute convergence of $\sum_{n=0}^{\infty}\left|c_{n}\right|$
is favored by large $\alpha$ and $N$, and is assured for $\alpha>1/2$
regardless of the value of $N$.

Where $\sum_{n=0}^{\infty}\left|c_{n}\right|$ converges, the Cauchy
product in Eq. (\ref{eq:p5}) not only converges but also does so
absolutely. Hence, by the dominated convergence theorem, we can take
expectations of $f^{2}$, i.e. $\boldsymbol{\mathbf{E}}f^{2}$ and
the corresponding approximations of $\boldsymbol{\mathbf{E}}f^{2}$
converge as well. Finally, since 
\begin{equation}
\sigma_{f}^{2}\boldsymbol{=\mathbf{E}}f^{2}-\left(\boldsymbol{\mathbf{E}}f\right)^{2}\leq\boldsymbol{=\mathbf{E}}f^{2}\label{eq:p13}
\end{equation}
the series expansion for $\sigma_{f}^{2}$ converges in case that
of $\boldsymbol{\mathbf{E}}f^{2}$ does. 

A sufficient condition for convergence of the series for variance
$\sigma_{f}^{2}$ is therefore $\frac{1-\alpha}{\alpha}<1$. In practice,
either owing to a small missing data probability so that $1-\alpha<1/2$,
and in its absence a large number of potential observations $N$ to
compensate, the series for variance converges. Therefore, we can use
practical moment-based estimators in geophysical missing-data problems.
The remainder of the section is devoted to deriving low-order estimators
for variance of the spatial mean. 

\subsection{Series for variance}

Here we shall derive low-order estimators for variance in the presence
of missing data, to study their structure. Truncating the series upto
$2$\textsuperscript{nd} order
\begin{multline}
f\left(R,S\right)=f\left(\boldsymbol{\mathbf{E}}R,\boldsymbol{\mathbf{E}}S\right)+\ldots+\frac{1}{2}\left(R-\boldsymbol{\mathbf{E}}R\right)^{2}\frac{\partial^{2}f}{\partial R^{2}}+\frac{1}{2}\left(S-\boldsymbol{\mathbf{E}}S\right)^{2}\frac{\partial^{2}f}{\partial S^{2}}+\\
\left(R-\boldsymbol{\mathbf{E}}R\right)\left(S-\boldsymbol{\mathbf{E}}S\right)\frac{\partial^{2}f}{\partial R\partial S},\label{eq:p14}
\end{multline}
with partial derivatives evaluated at $\left(\boldsymbol{\mathbf{E}}R,\boldsymbol{\mathbf{E}}S\right)$,
we obtain for the variance
\begin{multline}
\sigma_{f}^{2}=\boldsymbol{=\mathbf{E}}f^{2}-\left(\boldsymbol{\mathbf{E}}f\right)^{2}=-6\frac{\left(\boldsymbol{\mathbf{E}}R\right)^{2}}{\left(\boldsymbol{\mathbf{E}}S\right)^{2}}+\frac{4}{\left(\boldsymbol{\mathbf{E}}S\right)^{2}}\boldsymbol{\mathbf{E}}R^{2}+10\frac{\left(\boldsymbol{\mathbf{E}}R\right)^{2}}{\left(\boldsymbol{\mathbf{E}}S\right)^{4}}\boldsymbol{\mathbf{E}}S^{2}-\frac{1}{\left(\boldsymbol{\mathbf{E}}S\right)^{4}}\left(\boldsymbol{\mathbf{E}}RS\right)^{2}-\frac{\left(\boldsymbol{\mathbf{E}}R\right)^{2}}{\left(\boldsymbol{\mathbf{E}}S\right)^{6}}\left(\boldsymbol{\mathbf{E}}S^{2}\right)^{2}\\
+2\frac{\left(\boldsymbol{\mathbf{E}}R\right)}{\left(\boldsymbol{\mathbf{E}}S\right)^{4}}\boldsymbol{\mathbf{E}}RS^{2}-\frac{4}{\left(\boldsymbol{\mathbf{E}}S\right)^{3}}\boldsymbol{\mathbf{E}}R^{2}S-6\frac{\left(\boldsymbol{\mathbf{E}}R\right)^{2}}{\left(\boldsymbol{\mathbf{E}}S\right)^{5}}\boldsymbol{\mathbf{E}}S^{3}+\frac{\left(\boldsymbol{\mathbf{E}}R\right)^{2}}{\left(\boldsymbol{\mathbf{E}}S\right)^{6}}\boldsymbol{\mathbf{E}}S^{4}+\frac{1}{\left(\boldsymbol{\mathbf{E}}S\right)^{4}}\boldsymbol{\mathbf{E}}R^{2}S^{2}\\
+2\frac{\left(\boldsymbol{\mathbf{E}}R\right)\left(\boldsymbol{\mathbf{E}}RS\right)}{\left(\boldsymbol{\mathbf{E}}S\right)^{3}}\left\{ \frac{\boldsymbol{\mathbf{E}}S^{2}}{\left(\boldsymbol{\mathbf{E}}S\right)^{2}}-1\right\} .\label{eq:p15}
\end{multline}
These calculations are detailed in Supplementary Information Section
2. Evaluation of these statistics in general requires formulas for
$\boldsymbol{\mathbf{E}}R^{2}S^{l}$, $\boldsymbol{\mathbf{E}}RS^{l}$,
and $\boldsymbol{\mathbf{E}}S^{l}$, where $l$ is a non-negative
integer. 

\subsection{Evaluation of the moments}

These moments $\boldsymbol{\mathbf{E}}R^{2}S^{l}$, $\boldsymbol{\mathbf{E}}RS^{l}$,
and $\boldsymbol{\mathbf{E}}S^{l}$ $\left(l=0,1,2,\ldots\right)$
can be evaluated in a straightforward manner, but requires making
further assumptions. We assume that reporting of observations at different
locations is independent so that $\boldsymbol{\mathbf{E}}s_{i}s_{j}=\boldsymbol{\mathbf{E}}s_{i}\boldsymbol{\mathbf{E}}s_{j}$
if $i\neq j$. We also stipulate that availability is independent
of the measured field, so that $\boldsymbol{\mathbf{E}}r_{i}s_{i}=\boldsymbol{\mathbf{E}}r_{i}\boldsymbol{\mathbf{E}}s_{i}$.
Then, using the multinomial theorem, we obtain the moments as summarized
in Table 1. The derivations have been detailed in Supplementary Information
Section 3. 

From the multinomial theorem
\begin{equation}
\boldsymbol{\mathbf{E}}S^{l}=\boldsymbol{\mathbf{E}}\left(\sum_{i=1}^{N}\beta_{i}s_{i}\right)^{l}=\mathbf{E}\sum_{l_{1}+l_{2}+\ldots l_{N}=l}\frac{l!}{l_{1}!l_{2}!\cdots l_{N}!}\beta_{1}^{l_{1}}\beta_{2}^{l_{2}}\cdots\beta_{N}^{l_{N}}s_{1}^{l_{1}}s_{2}^{l_{2}}\cdots s_{N}^{l_{N}},
\end{equation}
and from linearity of expectation
\begin{equation}
\boldsymbol{\mathbf{E}}S^{l}=\sum_{l_{1}+l_{2}+\ldots l_{N}=l}\frac{l!}{l_{1}!l_{2}!\cdots l_{N}!}\beta_{1}^{l_{1}}\beta_{2}^{l_{2}}\cdots\beta_{N}^{l_{N}}\boldsymbol{\mathbf{E}}s_{1}^{l_{1}}s_{2}^{l_{2}}\cdots s_{N}^{l_{N}},
\end{equation}
and from independence of $s_{1}$, $s_{2}$, etc. 
\begin{equation}
\boldsymbol{\mathbf{E}}S^{l}=\sum_{l_{1}+l_{2}+\ldots l_{N}=l}\frac{l!}{l_{1}!l_{2}!\cdots l_{N}!}\beta_{1}^{l_{1}}\beta_{2}^{l_{2}}\cdots\beta_{N}^{l_{N}}\left(\boldsymbol{\mathbf{E}}s_{1}^{l_{1}}\boldsymbol{\mathbf{E}}s_{2}^{l_{2}}\cdots\mathbf{\boldsymbol{\mathbf{E}}}s_{N}^{l_{N}}\right).\label{eq:17v2}
\end{equation}

Now, $\boldsymbol{\mathbf{E}}s_{i}^{l_{i}}=\alpha$ if $l_{i}>0$
and $1$ otherwise. If there are $m$ distinct terms in the product
$\boldsymbol{\mathbf{E}}s_{1}^{l_{1}}\boldsymbol{\mathbf{E}}s_{2}^{l_{2}}\cdots\mathbf{\boldsymbol{\mathbf{E}}}s_{N}^{l_{N}}$
with nonzero power, then the product becomes $\alpha^{m}$. Similarly
\begin{equation}
\boldsymbol{\mathbf{E}}RS^{l}=\boldsymbol{\mathbf{E}}\sum_{i=1}^{N}\beta_{i}r_{i}s_{i}\sum_{l_{1}+l_{2}+\ldots l_{N}=l}\frac{l!}{l_{1}!l_{2}!\cdots l_{N}!}\beta_{1}^{l_{1}}\beta_{2}^{l_{2}}\cdots\beta_{N}^{l_{N}}s_{1}^{l_{1}}s_{2}^{l_{2}}\cdots s_{N}^{l_{N}}
\end{equation}
and from linearity of $\boldsymbol{\mathbf{E}}$ and independence
between $r_{i}$ and $s_{i}$
\begin{equation}
\boldsymbol{\mathbf{E}}RS^{l}=\sum_{i=1}^{N}\beta_{i}\boldsymbol{\mathbf{E}}r_{i}\left\{ \sum_{l_{1}+l_{2}+\ldots l_{N}=l}\frac{l!}{l_{1}!l_{2}!\cdots l_{N}!}\beta_{1}^{l_{1}}\beta_{2}^{l_{2}}\cdots\beta_{N}^{l_{N}}\left(\boldsymbol{\mathbf{E}}s_{1}^{l_{1}}\boldsymbol{\mathbf{E}}s_{2}^{l_{2}}\cdots\boldsymbol{\mathbf{E}}s_{N}^{l_{N}}\right)\boldsymbol{\mathbf{E}}s_{i}\right\} ,\label{eq:p55}
\end{equation}
where for fixed $i$ the product $\left(\boldsymbol{\mathbf{E}}s_{1}^{l_{1}}\boldsymbol{\mathbf{E}}s_{2}^{l_{2}}\cdots\boldsymbol{\mathbf{E}}s_{N}^{l_{N}}\right)\boldsymbol{\mathbf{E}}s_{i}$
equals $\alpha^{m}$ if there are $m$ distinct indices in the product
that include the $i$\textsuperscript{th} index. If one of the indices
among the nonzero $l_{k}$'s equals $i$ then there must be $m$ distinct
terms in the product $\boldsymbol{\mathbf{E}}s_{1}^{l_{1}}\boldsymbol{\mathbf{E}}s_{2}^{l_{2}}\cdots\boldsymbol{\mathbf{E}}s_{N}^{l_{N}}$,
otherwise there are $m-1$ distinct terms. 

Similarly we can compute $\boldsymbol{\mathbf{E}}R^{2}S^{l}$ as detailed
in Supplementary Information. For example, 
\begin{multline}
\boldsymbol{\mathbf{E}}R^{2}S^{2}=\alpha\sum\beta_{i}^{4}\boldsymbol{\mathbf{E}}r_{i}^{2}+\alpha^{2}\left\{ \sum_{j\neq i}\beta_{i}^{2}\beta_{j}^{2}\boldsymbol{\mathbf{E}}r_{i}^{2}+\sum_{j\neq i}2\beta_{i}^{3}\beta_{j}\boldsymbol{\mathbf{E}}r_{i}r_{j}+\sum_{j\neq i}2\beta_{i}^{3}\beta_{j}\boldsymbol{\mathbf{E}}r_{i}^{2}+\sum_{j\neq i}2\beta_{i}^{2}\beta_{j}^{2}\boldsymbol{\mathbf{E}}r_{i}r_{j}\right\} \\
+\alpha^{3}\left\{ \sum_{i\neq j\neq k}\beta_{i}\beta_{j}\beta_{k}^{2}\boldsymbol{\mathbf{E}}r_{i}r_{j}+\sum_{i\neq j\neq k}\beta_{i}^{2}\beta_{j}\beta_{k}\boldsymbol{\mathbf{E}}r_{i}^{2}+\sum_{i\neq j\neq k}4\beta_{i}^{2}\beta_{j}\beta_{k}\boldsymbol{\mathbf{E}}r_{i}r_{j}\right\} +\alpha^{4}\sum_{i\neq j\neq k\neq m}\beta_{i}\beta_{j}\beta_{k}\beta_{m}\boldsymbol{\mathbf{E}}r_{i}r_{j}.\label{eq:p16}
\end{multline}

Table 1: Moments appearing in truncated approximation for variance
in Eq. (\ref{eq:p15}).

\begin{tabular}{|c|c|}
\hline 
Moment & Expression\tabularnewline
\hline 
\hline 
$\boldsymbol{\mathbf{E}}S$ & $\alpha$\tabularnewline
\hline 
$\boldsymbol{\mathbf{E}}S^{2}$ & $\alpha\sum\beta_{k}^{2}+\alpha^{2}\sum_{j\neq k}\beta_{j}\beta_{k}$\tabularnewline
\hline 
$\boldsymbol{\mathbf{E}}S^{3}$ & $\alpha\sum\beta_{k}^{3}+\alpha^{2}\sum_{j\neq k}3\beta_{j}\beta_{k}^{2}+\alpha^{3}\sum_{i\neq j\neq k}\beta_{i}\beta_{j}\beta_{k}$\tabularnewline
\hline 
$\boldsymbol{\mathbf{E}}S^{4}$ & $\alpha\sum\beta_{k}^{4}+\alpha^{2}\left\{ \sum_{j\neq k}4\beta_{j}\beta_{k}^{3}+\sum_{j\neq k}3\beta_{j}^{2}\beta_{k}^{2}\right\} +\alpha^{3}\sum_{i\neq j\neq k}6\beta_{i}\beta_{j}\beta_{k}^{2}+\alpha^{4}\sum_{i\neq j\neq k\neq m}\beta_{i}\beta_{j}\beta_{k}\beta_{m}$\tabularnewline
\hline 
$\boldsymbol{\mathbf{E}}R$ & $\alpha\sum_{i=1}^{N}\beta_{i}\boldsymbol{\mathbf{E}}r_{i}$\tabularnewline
\hline 
$\boldsymbol{\mathbf{E}}RS$ & $\alpha\sum_{i=1}^{N}\beta_{i}^{2}\boldsymbol{\mathbf{E}}r_{i}+\alpha^{2}\sum_{i\neq j}\beta_{i}\beta_{j}\boldsymbol{\mathbf{E}}r_{i}$\tabularnewline
\hline 
$\boldsymbol{\mathbf{E}}RS^{2}$ & $\alpha\sum_{i=1}^{N}\beta_{i}^{3}\boldsymbol{\mathbf{E}}r_{i}+\alpha^{2}\left\{ \sum_{j\neq i}^{N}\beta_{i}\beta_{j}^{2}\boldsymbol{\mathbf{E}}r_{i}+\sum_{j\neq i}2\beta_{i}^{2}\beta_{j}\boldsymbol{\mathbf{E}}r_{i}\right\} +\alpha^{3}\sum_{i\neq j\neq k}2\beta_{i}\beta_{j}\beta_{k}\boldsymbol{\mathbf{E}}r_{i}$\tabularnewline
\hline 
$\boldsymbol{\mathbf{E}}R^{2}$ & $\alpha\sum\beta_{i}^{2}\boldsymbol{\mathbf{E}}r_{i}^{2}+\alpha^{2}\sum_{i\neq j}\beta_{i}\beta_{j}\boldsymbol{\mathbf{E}}r_{i}r_{j}$\tabularnewline
\hline 
$\boldsymbol{\mathbf{E}}R^{2}S$ & $\alpha\sum\beta_{i}^{3}\boldsymbol{\mathbf{E}}r_{i}^{2}+\alpha^{2}\left\{ \sum_{i\neq j}\beta_{i}^{2}\beta_{j}\boldsymbol{\mathbf{E}}r_{i}^{2}+2\sum_{i\neq j}\beta_{i}^{2}\beta_{j}\boldsymbol{\mathbf{E}}r_{i}r_{j}\right\} +\alpha^{3}\sum_{i\neq j\neq k}\beta_{i}\beta_{j}\beta_{k}\boldsymbol{\mathbf{E}}r_{i}r_{j}$\tabularnewline
\hline 
\end{tabular}

\section{Uniform weights}

The aforementioned formulas are applicable to general weights $\beta_{i}$
that sum to one. A simplification of wide importance results for uniform
weights $\beta_{i}=1/N$, as when arithmetic averages are taken. One
application of uniform weights is sampling without replacement to
make inferences about the variance of sample averages. This yields
simplified formulas for the moments. For example, from Table 1, for
uniform weights
\begin{equation}
\boldsymbol{\mathbf{E}}S^{2}=\alpha\sum\frac{1}{N^{2}}+\alpha^{2}\sum_{j\neq k}\frac{1}{N^{2}}=\alpha\frac{1}{N}+\alpha^{2}\frac{\left(N-1\right)}{N}\label{eq:p17}
\end{equation}
because the summations occur over $N$ indices and $N\left(N-1\right)$
indices respectively. Similarly, 
\begin{equation}
\boldsymbol{\mathbf{E}}R^{2}S=\alpha\sum\frac{1}{N^{3}}\boldsymbol{\mathbf{E}}r_{i}^{2}+\alpha^{2}\left\{ \sum_{i\neq j}\frac{1}{N^{3}}\boldsymbol{\mathbf{E}}r_{i}^{2}+\sum_{i\neq j}\frac{2}{N^{3}}\boldsymbol{\mathbf{E}}r_{i}r_{j}\right\} +\alpha^{3}\sum_{i\neq j\neq k}\frac{1}{N^{3}}\boldsymbol{\mathbf{E}}r_{i}r_{j},\label{eq:p18}
\end{equation}
which simplifies to 
\begin{equation}
\alpha\sum_{i}\frac{1}{N^{3}}\boldsymbol{\mathbf{E}}r_{i}^{2}+\alpha^{2}\left\{ \sum_{i}\frac{N-1}{N^{3}}\boldsymbol{\mathbf{E}}r_{i}^{2}+\sum_{i\neq j}\frac{2}{N^{3}}\boldsymbol{\mathbf{E}}r_{i}r_{j}\right\} +\alpha^{3}\sum_{i\neq j}\frac{N-2}{N^{3}}\boldsymbol{\mathbf{E}}r_{i}r_{j}.\label{eq:p19}
\end{equation}
Table 2 summarizes these moments. The simplifications are detailed
in Supplementary Information Section 4. As before, in addition to
$\alpha$, these depend only on the $2$\textsuperscript{nd} order
moments of the underlying field. 

\pagebreak{}

\begin{landscape}

Table 2: Moments for the case of uniform weights.

\begin{tabular}{|c|c|c|}
\hline 
Moment & Expression & Large $N$ approx.\tabularnewline
\hline 
\hline 
$\boldsymbol{\mathbf{E}}S$ & $\alpha$ & $\alpha$\tabularnewline
\hline 
$\boldsymbol{\mathbf{E}}S^{2}$ & $\alpha\frac{1}{N}+\alpha^{2}\frac{\left(N-1\right)}{N}$ & $\alpha^{2}$\tabularnewline
\hline 
$\boldsymbol{\mathbf{E}}S^{3}$ & $\boldsymbol{\mathbf{E}}S^{3}=\alpha\frac{1}{N^{2}}+3\alpha^{2}\frac{\left(N-1\right)}{N^{2}}+\alpha^{3}\frac{\left(N-1\right)\left(N-2\right)}{N^{2}}$ & $\alpha^{3}$\tabularnewline
\hline 
$\boldsymbol{\mathbf{E}}S^{4}$ & $\alpha\frac{1}{N^{3}}+7\alpha^{2}\frac{\left(N-1\right)}{N^{3}}+6\alpha^{3}\frac{\left(N-1\right)\left(N-2\right)}{N^{3}}+\alpha^{4}\frac{\left(N-1\right)\left(N-2\right)\left(N-3\right)}{N^{3}}$ & $\alpha^{4}$\tabularnewline
\hline 
$\boldsymbol{\mathbf{E}}R$ & $\alpha\frac{1}{N}\sum_{i=1}^{N}\boldsymbol{\mathbf{E}}r_{i}$ & $\alpha\frac{1}{N}\sum_{i=1}^{N}\boldsymbol{\mathbf{E}}r_{i}$\tabularnewline
\hline 
$\boldsymbol{\mathbf{E}}RS$ & $\left\{ \alpha\frac{1}{N^{2}}+\alpha^{2}\frac{N-1}{N^{2}}\right\} \sum_{i}\boldsymbol{\mathbf{E}}r_{i}$ & $\alpha^{2}\frac{1}{N}\sum_{i=1}^{N}\boldsymbol{\mathbf{E}}r_{i}$\tabularnewline
\hline 
$\boldsymbol{\mathbf{E}}RS^{2}$ & $\left\{ \alpha\frac{1}{N^{3}}+3\alpha^{2}\frac{N-1}{N^{3}}+\alpha^{3}\frac{\left(N-1\right)\left(N-2\right)}{N^{3}}\right\} \sum_{i}\boldsymbol{\mathbf{E}}r_{i}$ & $\alpha^{3}\frac{1}{N}\sum_{i=1}^{N}\boldsymbol{\mathbf{E}}r_{i}$\tabularnewline
\hline 
$\boldsymbol{\mathbf{E}}R^{2}$ & $\alpha\frac{1}{N^{2}}\sum_{i}\boldsymbol{\mathbf{E}}r_{i}^{2}+\alpha^{2}\frac{1}{N^{2}}\sum_{i\neq j}\boldsymbol{\mathbf{E}}r_{i}r_{j}$ & $\alpha\frac{1}{N^{2}}\sum_{i}\boldsymbol{\mathbf{E}}r_{i}^{2}+\alpha^{2}\frac{1}{N^{2}}\sum_{i\neq j}\boldsymbol{\mathbf{E}}r_{i}r_{j}$\tabularnewline
\hline 
$\boldsymbol{\mathbf{E}}R^{2}S$ & $\left\{ \alpha\frac{1}{N^{3}}+\alpha^{2}\frac{N-1}{N^{3}}\right\} \sum_{i}\boldsymbol{\mathbf{E}}r_{i}^{2}+\left\{ \alpha^{2}\frac{2}{N^{3}}+\alpha^{3}\frac{N-2}{N^{3}}\right\} \sum_{i\neq j}\boldsymbol{\mathbf{E}}r_{i}r_{j}$ & $\alpha^{2}\frac{1}{N^{2}}\sum_{i}\boldsymbol{\mathbf{E}}r_{i}^{2}+\alpha^{3}\frac{1}{N^{2}}\sum_{i\neq j}\boldsymbol{\mathbf{E}}r_{i}r_{j}$\tabularnewline
\hline 
$\boldsymbol{\mathbf{E}}R^{2}S^{2}$ & $\left\{ \alpha\frac{1}{N^{4}}+3\alpha^{2}\frac{N-1}{N^{4}}+\alpha^{3}\frac{\left(N-1\right)\left(N-2\right)}{N^{4}}\right\} \sum_{i}\boldsymbol{\mathbf{E}}r_{i}^{2}+\left\{ 4\alpha^{2}\frac{1}{N^{4}}+5\alpha^{3}\frac{N-2}{N^{4}}+\alpha^{4}\frac{\left(N-2\right)\left(N-3\right)}{N^{4}}\right\} \sum_{i\neq j}\boldsymbol{\mathbf{E}}r_{i}r_{j}$ & $\alpha^{3}\frac{1}{N^{2}}\sum_{i}\boldsymbol{\mathbf{E}}r_{i}^{2}+\alpha^{4}\frac{1}{N^{2}}\sum_{i\neq j}\boldsymbol{\mathbf{E}}r_{i}r_{j}$\tabularnewline
\hline 
\end{tabular}

\end{landscape}

\pagebreak{}

\subsection{Derivation using combinatorics}

These formulas for the case of uniform weights can be obtained more
simply using combinatorics. For uniform weights, Eq. (\ref{eq:17v2})
\begin{equation}
\boldsymbol{\mathbf{E}}S^{l}=\sum_{l_{1}+l_{2}+\ldots l_{N}=l}\frac{l!}{l_{1}!l_{2}!\cdots l_{N}!}\left(\frac{1}{N}\right)^{l_{1}}\left(\frac{1}{N}\right)^{l_{2}}\cdots\left(\frac{1}{N}\right)^{l_{N}}\left(\boldsymbol{\mathbf{E}}s_{1}^{l_{1}}\boldsymbol{\mathbf{E}}s_{2}^{l_{2}}\cdots\mathbf{\boldsymbol{\mathbf{E}}}s_{N}^{l_{N}}\right)\label{eq:p20}
\end{equation}
becomes 
\begin{equation}
\frac{1}{N^{l}}\sum_{l_{1}+l_{2}+\ldots l_{N}=l}\frac{l!}{l_{1}!l_{2}!\cdots l_{N}!}\left(\boldsymbol{\mathbf{E}}s_{1}^{l_{1}}\boldsymbol{\mathbf{E}}s_{2}^{l_{2}}\cdots\mathbf{\boldsymbol{\mathbf{E}}}s_{N}^{l_{N}}\right).\label{eq:p21}
\end{equation}
Consider the problem of assigning $l$ distinguishable balls in $N$
bins. This can be done in a total of $N^{l}$ different ways. This
is the sum of $\frac{l!}{l_{1}!l_{2}!\cdots l_{N}!}$, for $l_{1}+l_{2}+\ldots l_{N}=l$.
From the multinomial theorem, $\frac{l!}{l_{1}!l_{2}!\cdots l_{N}!}$
is the number of ways of maintaining $l_{1}$ balls in the $1$\textsuperscript{st}
bin, $l_{2}$ balls in the $2$\textsuperscript{nd} bin, etc., for
a total of $l$ balls in the $N$ bins. Now consider the problem of
assigning $l$ distinguishable balls in $N$ bins, such that exactly
$m$ bins are occupied, corresponding to $m$ distinct terms with
nonzero power in $\boldsymbol{\mathbf{E}}s_{1}^{l_{1}}\boldsymbol{\mathbf{E}}s_{2}^{l_{2}}\cdots\mathbf{\boldsymbol{\mathbf{E}}}s_{N}^{l_{N}}$. 

The number of ways of doing this corresponds to $N^{l}$ times the
coefficient of $\alpha^{m}$ in $\boldsymbol{\mathbf{E}}S^{l}$. This
coefficient is the number of ways of assigning $l$ distinguishable
balls in $N$ bins, such that exactly $m$ bins are occupied. This
is found in three steps: first, consider all the partitions of the
set of $l$ balls into exactly $m$ bins, with the number of such
partitions denoted by $\left\{ \begin{array}{c}
l\\
m
\end{array}\right\} $; second, choose the $m$ bins among the $N$ possible bins, for a
total of $\left(\begin{array}{c}
N\\
m
\end{array}\right)$ selections; and lastly, order the bins in one of $m!$ different
ways. Thus the coefficient describing the number of ways in which
$\boldsymbol{\mathbf{E}}s_{1}^{l_{1}}\boldsymbol{\mathbf{E}}s_{2}^{l_{2}}\cdots\mathbf{\boldsymbol{\mathbf{E}}}s_{N}^{l_{N}}$
can involve $m$ distinct terms, is given by
\begin{equation}
\left\{ \begin{array}{c}
l\\
m
\end{array}\right\} \left(\begin{array}{c}
N\\
m
\end{array}\right)m!=\left\{ \begin{array}{c}
l\\
m
\end{array}\right\} \frac{N!}{\left(N-m\right)!},\label{eq:p22}
\end{equation}
where $\left\{ \begin{array}{c}
l\\
m
\end{array}\right\} $, the number of partitions of $l$ distinguishable elements into exactly
$m$ nonempty sets, is given by Stirling numbers of the second kind
having formula
\begin{equation}
\left\{ \begin{array}{c}
l\\
m
\end{array}\right\} =\frac{1}{m!}\sum_{q=0}^{m}\left(-1\right)^{m-q}\left(\begin{array}{c}
m\\
q
\end{array}\right)q^{l}.\label{eq:p23}
\end{equation}
Finally we obtain
\begin{equation}
\boldsymbol{\mathbf{E}}S^{l}=\frac{1}{N^{l}}\sum_{m=1}^{\min\left\{ l,N\right\} }\left\{ \begin{array}{c}
l\\
m
\end{array}\right\} \frac{N!}{\left(N-m\right)!}\alpha^{m}\label{eq:p24}
\end{equation}
where the Stirling numbers $\left\{ \begin{array}{c}
l\\
m
\end{array}\right\} $ can be evaluated from the formula in Eq. (\ref{eq:p23}). Table 3
describes the coefficients of $\alpha^{m}$ in $\boldsymbol{\mathbf{E}}S^{l}$
for $1\leq l\leq4$, based on Eq. (\ref{eq:p24}), which confirms
the independently-derived formulas listed in Table 2. Corresponding
Stirling numbers of the second kind are documented in Table 7. 

\pagebreak{}

\begin{landscape}

Table 3: Coefficients of $\alpha^{m}$ in $\boldsymbol{\mathbf{E}}S^{l}$
, which is $\frac{1}{N^{l}}\left\{ \begin{array}{c}
l\\
m
\end{array}\right\} \frac{N!}{\left(N-m\right)!}$.

\begin{tabular}{|c|c|c|c|c|}
\hline 
 & $m=1$ & $m=2$ & $m=3$ & $m=4$\tabularnewline
\hline 
\hline 
$l=1$ & $\frac{1}{N^{1}}\left\{ \begin{array}{c}
1\\
1
\end{array}\right\} \frac{N!}{\left(N-1\right)!}=1$ &  &  & \tabularnewline
\hline 
$l=2$ & $\frac{1}{N^{2}}\left\{ \begin{array}{c}
2\\
1
\end{array}\right\} \frac{N!}{\left(N-1\right)!}=\frac{1}{N}$ & $\frac{1}{N^{2}}\left\{ \begin{array}{c}
2\\
2
\end{array}\right\} \frac{N!}{\left(N-2\right)!}=\frac{N-1}{N}$ &  & \tabularnewline
\hline 
$l=3$ & $\frac{1}{N^{3}}\left\{ \begin{array}{c}
3\\
1
\end{array}\right\} \frac{N!}{\left(N-1\right)!}=\frac{1}{N^{2}}$ & $\frac{1}{N^{3}}\left\{ \begin{array}{c}
3\\
2
\end{array}\right\} \frac{N!}{\left(N-2\right)!}=3\frac{N-1}{N^{2}}$ & $\frac{1}{N^{3}}\left\{ \begin{array}{c}
3\\
3
\end{array}\right\} \frac{N!}{\left(N-3\right)!}=\frac{\left(N-1\right)\left(N-2\right)}{N^{2}}$ & \tabularnewline
\hline 
$l=4$ & $\frac{1}{N^{4}}\left\{ \begin{array}{c}
4\\
1
\end{array}\right\} \frac{N!}{\left(N-1\right)!}=\frac{1}{N^{3}}$ & $\frac{1}{N^{4}}\left\{ \begin{array}{c}
4\\
2
\end{array}\right\} \frac{N!}{\left(N-2\right)!}=7\frac{N-1}{N^{3}}$ & $\frac{1}{N^{4}}\left\{ \begin{array}{c}
4\\
3
\end{array}\right\} \frac{N!}{\left(N-3\right)!}=6\frac{\left(N-1\right)\left(N-2\right)}{N^{3}}$ & $\frac{1}{N^{4}}\left\{ \begin{array}{c}
4\\
4
\end{array}\right\} \frac{N!}{\left(N-4\right)!}=\frac{\left(N-1\right)\left(N-2\right)\left(N-3\right)}{N^{3}}$\tabularnewline
\hline 
\end{tabular}

\end{landscape}

\pagebreak{}

Similarly, for uniform weights, Eq. (\ref{eq:p55}) 
\begin{equation}
\boldsymbol{\mathbf{E}}RS^{l}=\sum_{i=1}^{N}\frac{1}{N}\boldsymbol{\mathbf{E}}r_{i}\left\{ \sum_{l_{1}+l_{2}+\ldots l_{N}=l}\frac{l!}{l_{1}!l_{2}!\cdots l_{N}!}\left(\frac{1}{N}\right)^{l_{1}}\left(\frac{1}{N}\right)^{l_{2}}\cdots\left(\frac{1}{N}\right)^{l_{N}}\left(\boldsymbol{\mathbf{E}}s_{1}^{l_{1}}\boldsymbol{\mathbf{E}}s_{2}^{l_{2}}\cdots\boldsymbol{\mathbf{E}}s_{N}^{l_{N}}\right)\boldsymbol{\mathbf{E}}s_{i}\right\} \label{eq:p25}
\end{equation}
or
\begin{equation}
\boldsymbol{\mathbf{E}}RS^{l}=\frac{1}{N^{l+1}}\sum_{i=1}^{N}\boldsymbol{\mathbf{E}}r_{i}\left\{ \sum_{l_{1}+l_{2}+\ldots l_{N}=l}\frac{l!}{l_{1}!l_{2}!\cdots l_{N}!}\left(\boldsymbol{\mathbf{E}}s_{1}^{l_{1}}\boldsymbol{\mathbf{E}}s_{2}^{l_{2}}\cdots\boldsymbol{\mathbf{E}}s_{N}^{l_{N}}\right)\boldsymbol{\mathbf{E}}s_{i}\right\} ,\label{eq:p26}
\end{equation}

which is also a power-series in $\alpha$, where (by symmetry) $N$
times the coefficient of $\alpha^{m}$ is $\frac{1}{N^{l+1}}\sum_{i=1}^{N}\boldsymbol{\mathbf{E}}r_{i}$
times the number of ways of assigning $l+1$ distinguishable balls
in $N$ bins, such that exactly $m$ bins are occupied. The number
of ways of doing this is $\left\{ \begin{array}{c}
l+1\\
m
\end{array}\right\} \left(\begin{array}{c}
N\\
m
\end{array}\right)m!=\left\{ \begin{array}{c}
l+1\\
m
\end{array}\right\} \frac{N!}{\left(N-m\right)!}$, so that finally
\begin{equation}
\boldsymbol{\mathbf{E}}RS^{l}=\frac{1}{N^{l+1}}\sum_{i=1}^{N}\boldsymbol{\mathbf{E}}r_{i}\left\{ \sum_{m=1}^{\min\left\{ l+1,N\right\} }\frac{1}{N}\left\{ \begin{array}{c}
l+1\\
m
\end{array}\right\} \frac{N!}{\left(N-m\right)!}\alpha^{m}\right\} ,\label{eq:p27}
\end{equation}
which simplifies to
\begin{equation}
\boldsymbol{\mathbf{E}}RS^{l}=\frac{1}{N^{l+1}}\sum_{i=1}^{N}\boldsymbol{\mathbf{E}}r_{i}\left\{ \sum_{m=1}^{\min\left\{ l+1,N\right\} }\left\{ \begin{array}{c}
l+1\\
m
\end{array}\right\} \frac{\left(N-1\right)!}{\left(N-m\right)!}\alpha^{m}\right\} .\label{eq:p28-1}
\end{equation}
Coefficients of $\boldsymbol{\mathbf{E}}r_{i}\alpha^{m}$ in $\boldsymbol{\mathbf{E}}RS^{l}$
are shown in Table 4, which give the same results as corresponding
formulas in Table 2. 

\pagebreak{}

\begin{landscape}

Table 4: Coefficients of $\boldsymbol{\mathbf{E}}r_{i}\alpha^{m}$
in $\boldsymbol{\mathbf{E}}RS^{l}$, which is $\frac{1}{N^{l+1}}\left\{ \begin{array}{c}
l+1\\
m
\end{array}\right\} \frac{\left(N-1\right)!}{\left(N-m\right)!}$.

\begin{tabular}{|c|c|c|c|c|}
\hline 
 & $m=1$ & $m=2$ & $m=3$ & $m=4$\tabularnewline
\hline 
\hline 
$l=0$ & $\frac{1}{N^{1}}\left\{ \begin{array}{c}
1\\
1
\end{array}\right\} \frac{\left(N-1\right)!}{\left(N-1\right)!}=\frac{1}{N}$ &  &  & \tabularnewline
\hline 
$l=1$ & $\frac{1}{N^{2}}\left\{ \begin{array}{c}
2\\
1
\end{array}\right\} \frac{\left(N-1\right)!}{\left(N-1\right)!}=\frac{1}{N^{2}}$ & $\frac{1}{N^{2}}\left\{ \begin{array}{c}
2\\
2
\end{array}\right\} \frac{\left(N-1\right)!}{\left(N-2\right)!}=\frac{N-1}{N^{2}}$ &  & \tabularnewline
\hline 
$l=2$ & $\frac{1}{N^{3}}\left\{ \begin{array}{c}
3\\
1
\end{array}\right\} \frac{\left(N-1\right)!}{\left(N-1\right)!}=\frac{1}{N^{3}}$ & $\frac{1}{N^{3}}\left\{ \begin{array}{c}
3\\
2
\end{array}\right\} \frac{\left(N-1\right)!}{\left(N-2\right)!}=3\frac{N-1}{N^{3}}$ & $\frac{1}{N^{3}}\left\{ \begin{array}{c}
3\\
3
\end{array}\right\} \frac{\left(N-1\right)!}{\left(N-3\right)!}=\frac{\left(N-1\right)\left(N-2\right)}{N^{3}}$ & \tabularnewline
\hline 
$l=3$ & $\frac{1}{N^{4}}\left\{ \begin{array}{c}
4\\
1
\end{array}\right\} \frac{\left(N-1\right)!}{\left(N-1\right)!}=\frac{1}{N^{4}}$ & $\frac{1}{N^{4}}\left\{ \begin{array}{c}
4\\
2
\end{array}\right\} \frac{\left(N-1\right)!}{\left(N-2\right)!}=7\frac{N-1}{N^{4}}$ & $\frac{1}{N^{4}}\left\{ \begin{array}{c}
4\\
3
\end{array}\right\} \frac{\left(N-1\right)!}{\left(N-3\right)!}=6\frac{\left(N-1\right)\left(N-2\right)}{N^{4}}$ & $\frac{1}{N^{4}}\left\{ \begin{array}{c}
4\\
4
\end{array}\right\} \frac{\left(N-1\right)!}{\left(N-4\right)!}=\frac{\left(N-1\right)\left(N-2\right)\left(N-3\right)}{N^{4}}$\tabularnewline
\hline 
\end{tabular}

\end{landscape}

\pagebreak{}

Lastly $\boldsymbol{\mathbf{E}}R^{2}S^{l}$ becomes, for uniform weights
(from Eq. 23 in Supplementary Information)
\begin{multline}
\boldsymbol{\mathbf{E}}R^{2}S^{l}=\sum_{i=1}^{N}\frac{1}{N^{2}}\boldsymbol{\mathbf{E}}r_{i}^{2}\left\{ \sum_{l_{1}+l_{2}+\ldots l_{N}=l}\frac{l!}{l_{1}!l_{2}!\cdots l_{N}!}\left(\frac{1}{N}\right)^{l_{1}}\left(\frac{1}{N}\right)^{l_{2}}\cdots\left(\frac{1}{N}\right)^{l_{N}}\left(\boldsymbol{\mathbf{E}}s_{1}^{l_{1}}\boldsymbol{\mathbf{E}}s_{2}^{l_{2}}\cdots\boldsymbol{\mathbf{E}}s_{N}^{l_{N}}\right)\boldsymbol{\mathbf{E}}s_{i}^{2}\right\} \\
+\sum_{i\neq j}^{N}\frac{1}{N^{2}}\boldsymbol{\mathbf{E}}r_{i}r_{j}\left\{ \sum_{l_{1}+l_{2}+\ldots l_{N}=l}\frac{l!}{l_{1}!l_{2}!\cdots l_{N}!}\left(\frac{1}{N}\right)^{l_{1}}\left(\frac{1}{N}\right)^{l_{2}}\cdots\left(\frac{1}{N}\right)^{l_{N}}\left(\boldsymbol{\mathbf{E}}s_{1}^{l_{1}}\boldsymbol{\mathbf{E}}s_{2}^{l_{2}}\cdots\boldsymbol{\mathbf{E}}s_{N}^{l_{N}}\right)\boldsymbol{\mathbf{E}}s_{i}\boldsymbol{\mathbf{E}}s_{j}\right\} ,\label{eq:p29}
\end{multline}
or 
\begin{multline}
\boldsymbol{\mathbf{E}}R^{2}S^{l}=\frac{1}{N^{l+2}}\sum_{i=1}^{N}\boldsymbol{\mathbf{E}}r_{i}^{2}\left\{ \sum_{l_{1}+l_{2}+\ldots l_{N}=l}\frac{l!}{l_{1}!l_{2}!\cdots l_{N}!}\left(\boldsymbol{\mathbf{E}}s_{1}^{l_{1}}\boldsymbol{\mathbf{E}}s_{2}^{l_{2}}\cdots\boldsymbol{\mathbf{E}}s_{N}^{l_{N}}\right)\boldsymbol{\mathbf{E}}s_{i}^{2}\right\} \\
+\frac{1}{N^{l+2}}\sum_{i\neq j}^{N}\boldsymbol{\mathbf{E}}r_{i}r_{j}\left\{ \sum_{l_{1}+l_{2}+\ldots l_{N}=l}\frac{l!}{l_{1}!l_{2}!\cdots l_{N}!}\left(\boldsymbol{\mathbf{E}}s_{1}^{l_{1}}\boldsymbol{\mathbf{E}}s_{2}^{l_{2}}\cdots\boldsymbol{\mathbf{E}}s_{N}^{l_{N}}\right)\boldsymbol{\mathbf{E}}s_{i}\boldsymbol{\mathbf{E}}s_{j}\right\} ,\label{eq:p30}
\end{multline}
and since $\boldsymbol{\mathbf{E}}s_{i}^{2}=\boldsymbol{\mathbf{E}}s_{i}$
we can evaluate the first term directly as we did before. As for the
second term, we seek the number of distinct ways in which expressions
containing $\boldsymbol{\mathbf{E}}r_{i}r_{j}$ have $m$ distinct
terms in the product $\left(\boldsymbol{\mathbf{E}}s_{1}^{l_{1}}\boldsymbol{\mathbf{E}}s_{2}^{l_{2}}\cdots\boldsymbol{\mathbf{E}}s_{N}^{l_{N}}\right)\boldsymbol{\mathbf{E}}s_{i}\boldsymbol{\mathbf{E}}s_{j}$,
which corresponds to assignment of $l+2$ distinguishable balls in
$N$ bins, such that exactly $m$ bins are occupied. We must adjust
for the assignments already counted, involving terms of the form $\boldsymbol{\mathbf{E}}r_{i}^{2}$.
There are $2\left(\begin{array}{c}
N\\
2
\end{array}\right)$ ways of choosing each of the terms $\boldsymbol{\mathbf{E}}r_{i}r_{j}$,
$i\neq j$, and the terms in $r_{i}$'s are either of the form $\boldsymbol{\mathbf{E}}r_{i}^{2}$
or $\boldsymbol{\mathbf{E}}r_{i}r_{j}$. Therefore, denoting the coefficient
of $\alpha^{m}\boldsymbol{\mathbf{E}}r_{i}r_{j}$ as $b_{m}^{\left(l\right)}$and
that of $\alpha^{m}\boldsymbol{\mathbf{E}}r_{i}^{2}$ as $a_{m}^{\left(l\right)}$,
we obtain
\begin{equation}
Na_{m}^{\left(l\right)}+2\left(\begin{array}{c}
N\\
2
\end{array}\right)b_{m}^{\left(l\right)}=\frac{1}{N^{l+2}}\left\{ \begin{array}{c}
l+2\\
m
\end{array}\right\} \left(\begin{array}{c}
N\\
m
\end{array}\right)m!,\label{eq:p31}
\end{equation}
where, following the arguments above, the term on the right is the
number of ways to assign $l+2$ distinguishable balls into $m$ bins
chosen from among $N$ possible bins, and $Na_{m}^{\left(l\right)}=\frac{1}{N^{l+2}}\left\{ \begin{array}{c}
l+1\\
m
\end{array}\right\} \left(\begin{array}{c}
N\\
m
\end{array}\right)m!$. Therefore, $N\left(N-1\right)b_{m}^{\left(l\right)}=\frac{1}{N^{l+2}}\left(\left\{ \begin{array}{c}
l+2\\
m
\end{array}\right\} -\left\{ \begin{array}{c}
l+1\\
m
\end{array}\right\} \right)\left(\begin{array}{c}
N\\
m
\end{array}\right)m!$, and we obtain finally
\begin{multline}
\boldsymbol{\mathbf{E}}R^{2}S^{l}=\frac{1}{N^{l+2}}\sum_{i=1}\boldsymbol{\mathbf{E}}r_{i}^{2}\sum_{m=1}^{\min\left\{ l+1,N\right\} }\left\{ \begin{array}{c}
l+1\\
m
\end{array}\right\} \frac{\left(N-1\right)!}{\left(N-m\right)!}\alpha^{m}\\
+\frac{1}{N^{l+2}}\sum_{i\neq j}\boldsymbol{\mathbf{E}}r_{i}r_{j}\sum_{m=1}^{\min\left\{ l+2,N\right\} }\frac{\left(N-2\right)!}{\left(N-m\right)!}\left(\left\{ \begin{array}{c}
l+2\\
m
\end{array}\right\} -\left\{ \begin{array}{c}
l+1\\
m
\end{array}\right\} \right)\alpha^{m}.\label{eq:p32}
\end{multline}
These coefficients are shown in Tables 5 and 6, which confirm the
expressions listed in Table 2. In summary we have shown that, for
uniform weights, the moments appearing in the formula for variance
can be estimated via combinatorics. However, that does not give us
a formula for the variance itself. 

\pagebreak{}

\begin{landscape}

Table 5: Coefficients of $\boldsymbol{\mathbf{E}}r_{i}^{2}\alpha^{m}$
in $\boldsymbol{\mathbf{E}}R^{2}S^{l}$ , which is $\frac{1}{N^{l+2}}\left\{ \begin{array}{c}
l+1\\
m
\end{array}\right\} \frac{\left(N-1\right)!}{\left(N-m\right)!}$.

\begin{tabular}{|c|c|c|c|c|}
\hline 
 & $m=1$ & $m=2$ & $m=3$ & $m=4$\tabularnewline
\hline 
\hline 
$l=0$ & $\frac{1}{N^{2}}\left\{ \begin{array}{c}
1\\
1
\end{array}\right\} \frac{\left(N-1\right)!}{\left(N-1\right)!}=\frac{1}{N^{2}}$ &  &  & \tabularnewline
\hline 
$l=1$ & $\frac{1}{N^{3}}\left\{ \begin{array}{c}
2\\
1
\end{array}\right\} \frac{\left(N-1\right)!}{\left(N-1\right)!}=\frac{1}{N^{3}}$ & $\frac{1}{N^{3}}\left\{ \begin{array}{c}
2\\
2
\end{array}\right\} \frac{\left(N-1\right)!}{\left(N-2\right)!}=\frac{N-1}{N^{3}}$ &  & \tabularnewline
\hline 
$l=2$ & $\frac{1}{N^{4}}\left\{ \begin{array}{c}
3\\
1
\end{array}\right\} \frac{\left(N-1\right)!}{\left(N-1\right)!}=\frac{1}{N^{4}}$ & $\frac{1}{N^{4}}\left\{ \begin{array}{c}
3\\
2
\end{array}\right\} \frac{\left(N-1\right)!}{\left(N-2\right)!}=3\frac{N-1}{N^{4}}$ & $\frac{1}{N^{4}}\left\{ \begin{array}{c}
3\\
3
\end{array}\right\} \frac{\left(N-1\right)!}{\left(N-3\right)!}=\frac{\left(N-1\right)\left(N-2\right)}{N^{4}}$ & \tabularnewline
\hline 
$l=3$ & $\frac{1}{N^{5}}\left\{ \begin{array}{c}
4\\
1
\end{array}\right\} \frac{\left(N-1\right)!}{\left(N-1\right)!}=\frac{1}{N^{5}}$ & $\frac{1}{N^{5}}\left\{ \begin{array}{c}
4\\
2
\end{array}\right\} \frac{\left(N-1\right)!}{\left(N-2\right)!}=7\frac{N-1}{N^{5}}$ & $\frac{1}{N^{5}}\left\{ \begin{array}{c}
4\\
3
\end{array}\right\} \frac{\left(N-1\right)!}{\left(N-3\right)!}=6\frac{\left(N-1\right)\left(N-2\right)}{N^{5}}$ & $\frac{1}{N^{5}}\left\{ \begin{array}{c}
4\\
4
\end{array}\right\} \frac{\left(N-1\right)!}{\left(N-4\right)!}=\frac{\left(N-1\right)\left(N-2\right)\left(N-3\right)}{N^{5}}$\tabularnewline
\hline 
\end{tabular}

\pagebreak{}

Table 6: Coefficients of $\boldsymbol{\mathbf{E}}r_{i}r_{j}\alpha^{m}$
in $\boldsymbol{\mathbf{E}}R^{2}S^{l}$ , which is $\frac{1}{N^{l+2}}\frac{\left(N-2\right)!}{\left(N-m\right)!}\left(\left\{ \begin{array}{c}
l+2\\
m
\end{array}\right\} -\left\{ \begin{array}{c}
l+1\\
m
\end{array}\right\} \right)$.

\begin{tabular}{|c|c|c|c|c|}
\hline 
 & $m=1$ & $m=2$ & $m=3$ & $m=4$\tabularnewline
\hline 
\hline 
$l=0$ & $\frac{1}{N^{2}}\frac{\left(N-2\right)!}{\left(N-1\right)!}\left[\left\{ \begin{array}{c}
2\\
1
\end{array}\right\} -\left\{ \begin{array}{c}
1\\
1
\end{array}\right\} \right]=0$ & $\frac{1}{N^{2}}\frac{\left(N-2\right)!}{\left(N-2\right)!}\left[\left\{ \begin{array}{c}
2\\
2
\end{array}\right\} -\left\{ \begin{array}{c}
1\\
2
\end{array}\right\} \right]=\frac{1}{N^{2}}$ &  & \tabularnewline
\hline 
$l=1$ & $\frac{1}{N^{3}}\frac{\left(N-2\right)!}{\left(N-1\right)!}\left[\left\{ \begin{array}{c}
3\\
1
\end{array}\right\} -\left\{ \begin{array}{c}
2\\
1
\end{array}\right\} \right]=0$ & $\frac{1}{N^{3}}\frac{\left(N-2\right)!}{\left(N-2\right)!}\left[\left\{ \begin{array}{c}
3\\
2
\end{array}\right\} -\left\{ \begin{array}{c}
2\\
2
\end{array}\right\} \right]=\frac{2}{N^{3}}$ & $\frac{1}{N^{3}}\frac{\left(N-2\right)!}{\left(N-3\right)!}\left[\left\{ \begin{array}{c}
3\\
3
\end{array}\right\} -\left\{ \begin{array}{c}
2\\
3
\end{array}\right\} \right]=\frac{N-2}{N^{3}}$ & \tabularnewline
\hline 
$l=2$ & $\frac{1}{N^{4}}\frac{\left(N-2\right)!}{\left(N-1\right)!}\left[\left\{ \begin{array}{c}
4\\
1
\end{array}\right\} -\left\{ \begin{array}{c}
3\\
1
\end{array}\right\} \right]=0$ & $\frac{1}{N^{4}}\frac{\left(N-2\right)!}{\left(N-2\right)!}\left[\left\{ \begin{array}{c}
4\\
2
\end{array}\right\} -\left\{ \begin{array}{c}
3\\
2
\end{array}\right\} \right]=\frac{4}{N^{4}}$ & $\frac{1}{N^{4}}\frac{\left(N-2\right)!}{\left(N-3\right)!}\left[\left\{ \begin{array}{c}
4\\
3
\end{array}\right\} -\left\{ \begin{array}{c}
3\\
3
\end{array}\right\} \right]=5\frac{N-2}{N^{4}}$ & $\frac{1}{N^{4}}\frac{\left(N-2\right)!}{\left(N-4\right)!}\left[\left\{ \begin{array}{c}
4\\
4
\end{array}\right\} -\left\{ \begin{array}{c}
3\\
4
\end{array}\right\} \right]=\frac{\left(N-2\right)\left(N-3\right)}{N^{4}}$\tabularnewline
\hline 
\end{tabular}

\pagebreak{}

\end{landscape}

Table 7: Stirling numbers of the second kind $\left\{ \begin{array}{c}
n\\
k
\end{array}\right\} $, describing the number of partitions of $n$ elements into $k$ sets,
for the ranges used in Tables 3-6. 

\begin{tabular}{|c|c|c|c|c|c|}
\hline 
 & $k$ & $1$ & $2$ & $3$ & $4$\tabularnewline
\hline 
\hline 
$n$ &  &  &  &  & \tabularnewline
\hline 
$1$ &  & $1$ & $0$ & $0$ & $0$\tabularnewline
\hline 
$2$ &  & $1$ & $1$ & $0$ & $0$\tabularnewline
\hline 
$3$ &  & $1$ & $3$ & $1$ & $0$\tabularnewline
\hline 
$4$ &  & $1$ & $7$ & $6$ & $1$\tabularnewline
\hline 
\end{tabular}

\subsection{Large $N$ approximations}

These formulas become simplified if the field is sampled from a very
large collection of potential sites because, for fixed $\alpha$,
we need only consider the term with highest power of $N$ in the numerator.
For example, in the limit that $N\rightarrow\infty$, $\boldsymbol{\mathbf{E}}S^{4}\approx\alpha^{4}\frac{\left(N-1\right)\left(N-2\right)\left(N-3\right)}{N^{3}}\approx\alpha^{4}=\boldsymbol{\mathbf{E}}S^{4}.$
Similarly $\boldsymbol{\mathbf{E}}R^{2}S^{2}\approx\frac{\alpha^{3}}{N^{2}}\sum_{i}\boldsymbol{\mathbf{E}}r_{i}^{2}+\frac{\alpha^{4}}{N^{2}}\sum_{j\neq i}\boldsymbol{\mathbf{E}}r_{i}r_{j}=\boldsymbol{\mathbf{E}}R^{2}\boldsymbol{\mathbf{E}}S^{2}$,
and so on. Each of these simplifications is listed alongside the complete
formulas in Table 2. In the large $N$ limit, assuming $l\ll N$ we
obtain 
\begin{equation}
\boldsymbol{\mathbf{E}}S^{l}\approx\frac{1}{N^{l}}\left\{ \begin{array}{c}
l\\
l
\end{array}\right\} \frac{N!}{\left(N-l\right)!}\alpha^{l}=\frac{\left(N-1\right)\ldots\left(N-l+1\right)}{N^{l}}\alpha^{l}\approx\alpha^{l}=\left(\boldsymbol{\mathbf{E}}S\right)^{l}.\label{eq:p33}
\end{equation}
Similarly
\begin{equation}
\boldsymbol{\mathbf{E}}RS^{l}\approx\frac{1}{N^{l+1}}\sum_{i=1}^{N}\boldsymbol{\mathbf{E}}r_{i}\left\{ \begin{array}{c}
l+1\\
l+1
\end{array}\right\} \frac{\left(N-1\right)!}{\left(N-l-1\right)!}\alpha^{l+1}\approx\frac{1}{N}\sum_{i=1}^{N}\boldsymbol{\mathbf{E}}r_{i}\alpha^{l+1}=\boldsymbol{\mathbf{E}}R\left(\boldsymbol{\mathbf{E}}S\right)^{l},\label{eq:p34}
\end{equation}
whereas
\begin{multline}
\boldsymbol{\mathbf{E}}R^{2}S^{l}\approx\frac{1}{N^{l+2}}\sum_{i=1}\boldsymbol{\mathbf{E}}r_{i}^{2}\left\{ \begin{array}{c}
l+1\\
l+1
\end{array}\right\} \frac{\left(N-1\right)!}{\left(N-l-1\right)!}\alpha^{l+1}+\frac{1}{N^{l+2}}\sum_{i\neq j}\boldsymbol{\mathbf{E}}r_{i}r_{j}\frac{\left(N-2\right)!}{\left(N-l-2\right)!}\left(\left\{ \begin{array}{c}
l+2\\
l+2
\end{array}\right\} -\left\{ \begin{array}{c}
l+1\\
l+2
\end{array}\right\} \right)\alpha^{l+2}\\
\approx\left\{ \frac{\alpha}{N^{2}}\sum_{i=1}\boldsymbol{\mathbf{E}}r_{i}^{2}+\frac{\alpha^{2}}{N^{2}}\sum_{i\neq j}\boldsymbol{\mathbf{E}}r_{i}r_{j}\right\} \alpha^{l}=\boldsymbol{\mathbf{E}}R^{2}\boldsymbol{\mathbf{E}}S^{l}.\label{eq:p35}
\end{multline}

Substituting these results into Eq. (\ref{eq:p15}), the formula for
variance becomes
\begin{equation}
\sigma_{f}^{2}=\frac{1}{\alpha^{2}}\left\{ \boldsymbol{\mathbf{E}}R^{2}-\left(\boldsymbol{\mathbf{E}}R\right)^{2}\right\} =\frac{1}{\alpha^{2}}\sigma_{R}^{2},\label{eq:p36}
\end{equation}
which also equals the 1st order approximation Eq. (\ref{eq:p12})
to the variance when the moments of $R$ and $S$ are statistically
independent. The variance is of the ratio in Eq. (\ref{eq:p1}) is
simply the variance of the numerator inflated by $1/\alpha^{2}$ owing
to missing observations. 

\subsection{Behavior for moderate values of $N$ }

The behavior is far more interesting if $N$ is not necessarily large.
If we relax our assumption of $N\rightarrow\infty$, the moments of
$R$ and $S$ covary, and we must refer to the middle column of Table
2, for the case of $\beta_{i}=1/N$. Then the $2$\textsuperscript{nd}
order approximation for variance simplifies to
\begin{multline}
\sigma_{f}^{2}=\frac{1}{\alpha}\left[1+\frac{1}{N}\frac{\alpha-1}{\alpha}+\frac{1}{N^{2}}\frac{\alpha-1}{\alpha}\frac{2\alpha-1}{\alpha}\right]\frac{1}{N^{2}}\sum_{i}\boldsymbol{\mathbf{E}}r_{i}^{2}\\
+\left[1+\frac{3}{N}\frac{\alpha-1}{\alpha}+\frac{1}{N^{2}}\frac{\alpha-1}{\alpha}\frac{6\alpha-4}{\alpha}\right]\frac{1}{N^{2}}\sum_{j\neq i}\boldsymbol{\mathbf{E}}r_{i}r_{j}\\
-\left[1+\frac{2}{N}\frac{\alpha-1}{\alpha}-\frac{1}{N^{2}}\left(\frac{\alpha-1}{\alpha}\right)\left(\frac{3\alpha-1}{\alpha}\right)+\frac{1}{N^{3}}\left\{ 6\left(\frac{\alpha-1}{\alpha}\right)^{2}+\frac{1}{\alpha^{2}}\frac{\alpha-1}{\alpha}\right\} \right]\left\{ \frac{1}{N}\sum_{i=1}^{N}\boldsymbol{\mathbf{E}}r_{i}\right\} ^{2}.\label{eq:p38}
\end{multline}
Clearly each term is smaller by a factor of $\frac{1}{N}\left|\frac{\alpha-1}{\alpha}\right|$.
Recall from Section 2.1 that the inverse of this factor, i.e. $N\left|\frac{\alpha}{1-\alpha}\right|$
is the squared number of standard deviations from its expectation
$\boldsymbol{\mathbf{E}}S$ that $S$ must be, for the series to not
converge. The larger is this quantity, the more rapidly convergence
occurs. Therefore, we can expect rapid convergence of the series when
$\alpha$ is close to $1$, and more generally when $N$ is large
to offset the effect of $\left|\frac{\alpha-1}{\alpha}\right|>1$
when $\alpha$ is small. 

\subsection*{Simplifications}

Let us consider the following simplified cases of the variance formula
in Eq. (\ref{eq:p38}):
\begin{itemize}
\item $\alpha=1$: In case observations are reported with certainty, then
the variance 
\begin{equation}
\sigma_{f}^{2}=\frac{1}{N^{2}}\sum_{i}\boldsymbol{\mathbf{E}}r_{i}^{2}+\frac{1}{N^{2}}\sum_{j\neq i}\boldsymbol{\mathbf{E}}r_{i}r_{j}-\left\{ \frac{1}{N}\sum_{i=1}^{N}\boldsymbol{\mathbf{E}}r_{i}\right\} ^{2}
\end{equation}
simplifies to
\begin{equation}
\sigma_{f}^{2}=\frac{1}{N^{2}}\sum_{i}\sigma_{i}^{2}+\frac{1}{N^{2}}\sum_{j\neq i}\sigma_{ij},
\end{equation}
where $\sigma_{i}^{2}=\boldsymbol{\mathbf{E}}r_{i}^{2}-\left(\boldsymbol{\mathbf{E}}r_{i}\right)^{2}$
and $\sigma_{ij}=\boldsymbol{\mathbf{E}}r_{i}r_{j}-\boldsymbol{\mathbf{E}}r_{i}\boldsymbol{\mathbf{E}}r_{j}$
are the point-wise variance and pair-wise covariance. The formula
can be written equivalently as $\sigma_{i}^{2}=\beta^{T}\Sigma\beta$,
where $\Sigma$ is the covariance matrix, and $\beta=\left\{ \begin{array}{ccc}
\beta_{1} & \ldots & \beta_{n}\end{array}\right\} ^{T}$ is the vector of weights equal to $\left\{ \begin{array}{ccc}
\frac{1}{N} & \ldots & \frac{1}{N}\end{array}\right\} ^{T}$. 
\item $\alpha\rightarrow1$: When the probability of reporting is close
to $1$, we may consider only the first term in each of the series
\begin{equation}
\sigma_{f}^{2}=\frac{1}{\alpha}\frac{1}{N^{2}}\sum_{i}\boldsymbol{\mathbf{E}}r_{i}^{2}+\frac{1}{N^{2}}\sum_{j\neq i}\boldsymbol{\mathbf{E}}r_{i}r_{j}-\left\{ \frac{1}{N}\sum_{i=1}^{N}\boldsymbol{\mathbf{E}}r_{i}\right\} ^{2}
\end{equation}
which becomes
\begin{equation}
\sigma_{f}^{2}=\frac{1}{\alpha}\frac{1}{N^{2}}\sum_{i}\sigma_{i}^{2}+\frac{1}{N^{2}}\sum_{j\neq i}\sigma_{ij}+\frac{1-\alpha}{\alpha}\frac{1}{N}\sum_{i=1}^{N}\left(\boldsymbol{\mathbf{E}}r_{i}\right)^{2}.
\end{equation}
Compared to the previous case, the point-wise variance is inflated
by $1/\alpha$ and there is an additional correction proportional
to $\frac{1-\alpha}{\alpha}$ for squared point-wise expectations. 
\item Single epoch: The formula can be applied to a single epoch to describe
sampling without replacement of potential observations at a given
time. Expectations reduce to the variables themselves, so that $\boldsymbol{\mathbf{E}}r_{i}^{2}=r_{i}^{2},$
$\boldsymbol{\mathbf{E}}r_{i}r_{j}=r_{i}r_{j}$, and $\boldsymbol{\mathbf{E}}r_{i}=r_{i}$.
The formula for variance becomes
\begin{multline}
\sigma_{f}^{2}=\frac{1-\alpha}{\alpha}\frac{1}{N^{2}}\left[1+\frac{1}{N}\frac{2\alpha-1}{\alpha}-\frac{1}{N^{2}}\frac{\alpha^{2}+\alpha-1}{\alpha^{2}}+\frac{1}{N^{3}}\frac{6\alpha^{2}-6\alpha+1}{\alpha^{2}}\right]\sum_{i}r_{i}^{2}\\
-\frac{1-\alpha}{\alpha}\frac{1}{N^{3}}\left[1+\frac{1}{N}\frac{7\alpha-5}{\alpha}-\frac{1}{N^{2}}\frac{6\alpha^{2}-6\alpha+1}{\alpha^{2}}\right]\sum_{j\neq i}r_{i}r_{j}\label{eq:p43}
\end{multline}
In case $N$ is large so that both $\frac{1}{N}\left|\frac{2\alpha-1}{\alpha}\right|\ll1$
and $\frac{1}{N}\left|\frac{7\alpha-5}{\alpha}\right|\ll1$ , we can
simplify the above formula
\begin{equation}
\sigma_{f}^{2}=\frac{1-\alpha}{\alpha}\frac{1}{N^{2}}\left\{ \sum_{i}r_{i}^{2}-\frac{1}{N}\sum_{j\neq i}r_{i}r_{j}\right\} 
\end{equation}
and, for large $N$, $\sum_{j\neq i}r_{i}r_{j}\cong\left(\sum_{i}r_{i}\right)^{2}$,
so defining spatial variance for the epoch
\begin{equation}
\sigma_{s}^{2}=\frac{1}{N}\sum_{i}r_{i}^{2}-\left(\frac{1}{N}\sum_{i}r_{i}\right)^{2}
\end{equation}
we get 
\begin{equation}
\sigma_{f}^{2}=\frac{1-\alpha}{\alpha}\frac{1}{N}\sigma_{s}^{2}\label{eq:p50}
\end{equation}
as derived previously (\citet{Seshadri2018,Prakash2019}). Thus, for
a large number of observations, the ensemble variance in case of sampling
without replacement is proportional to the variance among the potential
samples. However, this simple result is merely a special case and
more generally one requires the formula in Eq. (\ref{eq:p43}).
\end{itemize}

\section{Simulations}

Illustrative simulations of a gridded rainfall dataset produced by
the Indian Meteorological Department have been undertaken, for a single
epoch. The dataset has daily rainfall over 357 gridpoints over India,
and we average this for one year (2001) over the summer monsoon months,
for 357 spatial values. We select a random subset of these gridded
values (of size $N$ ranging from $10-300$), fixing $N.$ For each
$N,$ an ensemble of Bernoulli random variables $\left\{ s_{i}:1\leq i\leq N\right\} $
is simulated. In each ensemble, each $s_{i}$ has probability $\alpha$
of equaling $1$. This indicates which of these gridded values are
employed in the spatial average calculation. Each ensemble member
has an associated spatial average, and the ensemble variance is computed,
and compared with the derived formula in Eq. (\ref{eq:p50}) (Figure
1). Each point describes one combination of $N$ and $\alpha.$ It
is seen in the left panel that the simulated variance departs from
the formula when the variance is large, i.e. $\frac{1-\alpha}{\alpha}\frac{1}{N}$
is large. The right panel shows the magnitude of the relative error,
and contour line where this takes the value $0.1,$ which distinguishes
regimes where the approximation is adequate/poor. Small $N$ is adequate
for $\alpha>1/2$ whereas much larger datasets are needed for small
$\alpha$. 

In practice $\frac{1-\alpha}{\alpha}\frac{1}{N}$ alone is not enough
to know the rates of convergence, but this factor recurs in the various
terms in the expansion. Figure 2 shows the two correction terms (upto
$1/N^{2}$) as functions of $\alpha$ and $N$, appearing in Eq. (\ref{eq:p38}).
Since these are corrections to the coefficient of $1$, they can be
neglected only for large $\alpha$ or $N.$ The second correction
term is generally negligible, except for small $\alpha$. Therefore
in case the missing data probability $1-\alpha$ is small one can
expect low-order truncations of the variance estimator to work quite
well. 

\pagebreak{}

\begin{figure}
\includegraphics[scale=0.55]{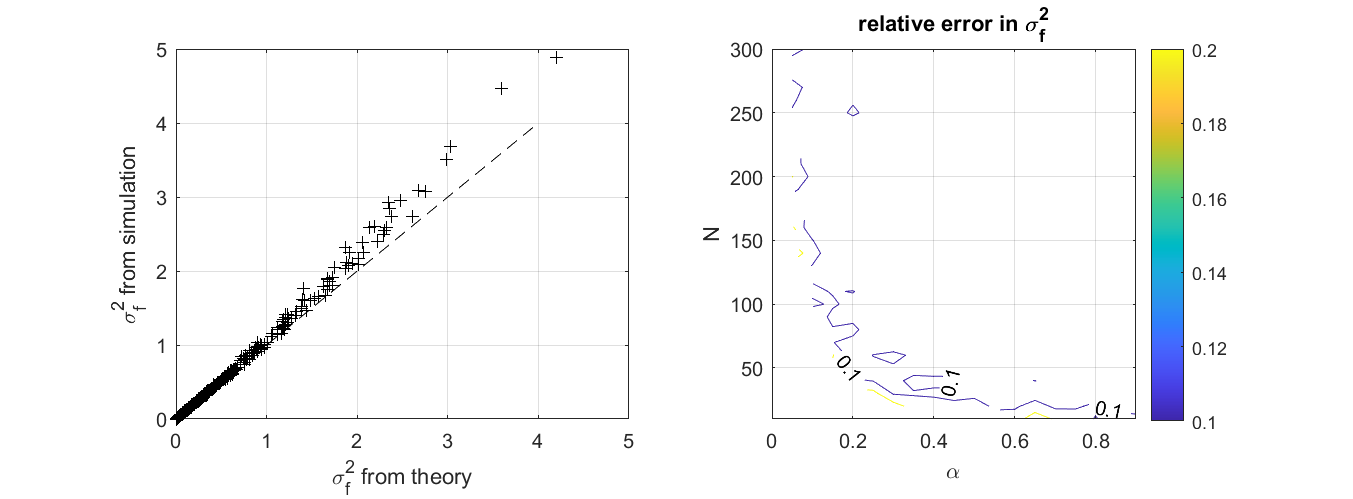}

\caption{Variance for a single epoch from Monte-Carlo simulations, compared
with the derived formula in Eq. (\ref{eq:p50}), in the left panel,
and the relative error (absolute value) in the right panel. The contour
where the relative error is $0.1$ is indicated. }
\end{figure}

\pagebreak{}

\begin{figure}
\includegraphics[scale=0.5]{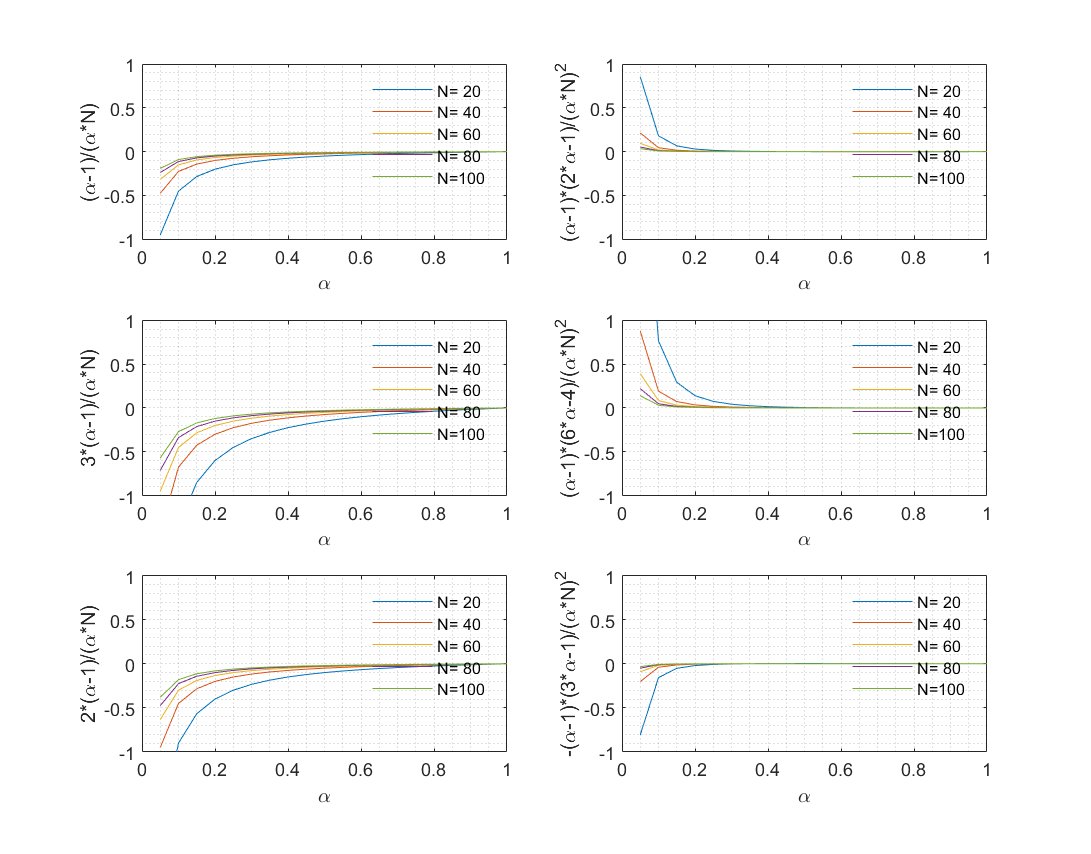}\caption{First two correction terms in the coefficients of $\sum_{i}\boldsymbol{\mathbf{E}}r_{i}^{2}$
(top panels), $\sum_{j\protect\neq i}\boldsymbol{\mathbf{E}}r_{i}r_{j}$
(middle panels), and $\left\{ \frac{1}{N}\sum_{i=1}^{N}\boldsymbol{\mathbf{E}}r_{i}\right\} ^{2}$
(bottom panels), as functions of $\alpha$ and $N$, from Eq. (\ref{eq:p38}).
Since these are corrections to the coefficient of $1$, they are generally
small for large $\alpha$ or $N.$}

\end{figure}

\pagebreak{}

\section{Discussion}

Deriving convergent estimators for the variance of a spatial mean
has broad relevance in the geosciences. These variance estimators
depend only on the first and second moments of the physical field,
and their convergence requires these to be finite. Many geophysical
problems experience finite first and second moments, so convergent
variance estimators can be derived for a broad range of settings,
from hydrology to climate science. Convergence is furthermore assured
if the probability of reporting individual observations is sufficiently
large (i.e., $\alpha>1/2)$. We have illustrated this convergence
by describing the behavior of successive terms in the series expansion
for variance. The results indicate how we can get good estimators
of the variance even with small datasets, as long as the probability
of missing data is not too large. Alternately, if the chance of reporting
an observation is small, the dataset had better be large (i.e., large
$N$).

An important special case occurs if the weights are all uniform, as
when arithmetic averages of the population are sought. This is sampling
without replacement from the population of potential observation sites
to make inferences about the properties of spatial averages. As in
the present study, if the probability of individual sites being observed
is held fixed, the total number of reported observations varies and
follows the Binomial distribution. While this differs from the Bootstrap
(\citet{Diaconis1983}) in its design and in sampling without replacement,
it still gives rise to useful nonparametric formulas of estimators
that arise in analogous contexts (\citet{Hankin2019}), and involving
statistics of subsamples from a population. Since the underlying theory
is nonparametric and holds general lessons, further investigations
of this situation is useful even though the exact formulation is different
from the more typical treatments. 

In the second part of this paper, we undertook a detailed analysis
of this situation. Compared to the general problem, this special case
with uniform weights inherits additional structure from our assumption
of the independence of the $s_{i}$'s. As a consequence the moments
in this case can be estimated directly through combinatorics. In particular,
all coefficients of the moments appearing in the estimators can be
easily computed via Stirling numbers of the second kind that count
the number of distinct set partitions (\citet{Aigner2007}). The problem
of convergence of the approximate estimators is also made explicit
as a result. Extending the analysis to higher approximations of the
variance might be useful. 

\bibliographystyle{agufull08}
\bibliography{Spatial_average_refs}

\end{document}